\documentclass[12pt]{article}
\usepackage{textcomp}
\usepackage{graphicx}           % 插入图片
\usepackage{epstopdf}
\usepackage{color}              % 支持彩色
\usepackage{xcolor,graphicx}
\usepackage{geometry}
\usepackage{float}
\usepackage{indentfirst}        % 首行缩进宏包
\usepackage{amsmath,amssymb,bm} % AMS math 宏包与数学符号加粗
\usepackage{amsthm}
\usepackage{cases}              % 数学公式
\usepackage{setspace}
\usepackage{hyperref}
\usepackage{titlesec}
\usepackage[all]{xy}            % 交换图
\usepackage{tikz}
\usetikzlibrary{arrows.meta}
\usepackage{multirow}% 合并 Table 的单元格时用到
\usepgflibrary{arrows}
\hypersetup{dvips,ps2pdf,CJKbookmarks,hypertex,
	colorlinks=true,
	pdfpagemode=FullScreen,
	bookmarksopen=true,
	bookmarksnumbered=true}     % 在中文环境下设置超链接文本
\usepackage{graphicx} %use graph format
\usepackage{epstopdf}
\usepackage{caption}
\usepackage{diagbox}
\usepackage{mathrsfs}
\usepackage{tikz}
\usepackage{dsfont}
\usepackage{graphics}
\usepackage{amsmath,amssymb,amsfonts,amsthm,bm,mathrsfs}
\usepackage{epsfig,graphicx,epstopdf,picinpar,subfigure,pstricks}
\usepackage{hyperref,cleveref,cite}
\usepackage{titlesec}
\usepackage{tcolorbox}
\usepackage{fancyhdr,fancyvrb}
\usepackage{xcolor,graphicx,geometry,float,indentfirst,cases,setspace,titlesec}
\usepackage{diagbox}

\newtheorem{thm}{Theorem}[section]

\newtheorem{coro}[thm]{Corollary}
\newtheorem{con}[thm]{Conjecture}

\newtheorem{pf}{proof}[section]
\theoremstyle{remark}

\begingroup
\makeatletter \@addtoreset{equation}{section} \makeatother
\makeindex 
\def\qed{\hfill \rule{4pt}{7pt}}

\begin{document}

\begin{center}

 {\Large \bf  Tur\'an inequalities for the broken $k$-diamond \\[5pt]
 partition function}

\end{center}

\begin{center}
Janet J.W. Dong$^{1}$, {Kathy Q. Ji}$^{2}$ and {Dennis X.Q. Jia}$^{3}$  \vskip 2mm

$^{1,\,2,\,3}$ Center for Applied Mathematics\\[3pt]
Tianjin University\\[3pt]
Tianjin 300072, P.R. China\\[6pt]
   \vskip 2mm

 Emails: $^1$dongjinwei@tju.edu.cn,   $^2$kathyji@tju.edu.cn and $^3$jxqmail@tju.edu.cn
\end{center}

\vskip 6mm \noindent {\bf Abstract.} We obtain an asymptotic formula for Andrews and Paule's   broken $k$-diamond partition function $\Delta_k(n)$ where $k=1$ or $2$. Based on this asymptotic formula, we  derive that  $\Delta_k(n)$ satisfies the order  $d$ Tur\'an inequalities for $d\geq 1$ and for sufficiently large $n$ when $k=1$ and $ 2$ by using a general result of Griffin, Ono, Rolen and Zagier. We also   show that Andrews and Paule's   broken $k$-diamond partition function $\Delta_k(n)$ is log-concave for $n\geq 1$ when $k=1$ and $2$.   This  leads to  $\Delta_k(a)\Delta_k(b)\ge\Delta_k(a+b)$ for $a,b\ge 1$ when $k=1$ and $ 2$.



 \vskip 6mm

\section{Introduction}

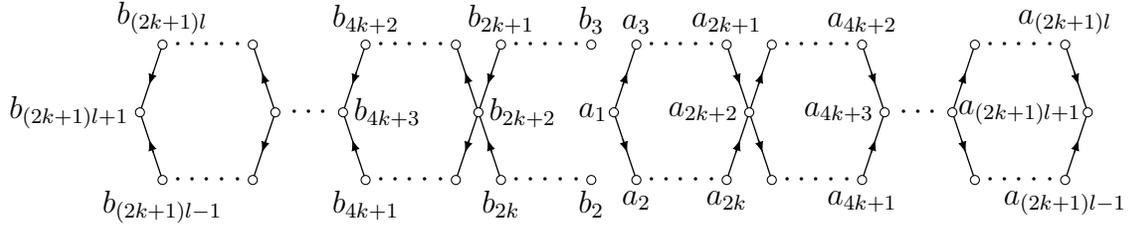
\begin{figure}[htbp]
	\begin{center}
		\begin{tikzpicture}[>=Stealth,scale=0.6]	
			\draw (-10.5,0) circle(.1);
			\node at (-10.5,0) [left]{$b_{(2k+1)l+1}$};
			\draw (0,0) circle(.1);
			\node at (0.1,0) [left] {$a_1$};
			\draw (3,0) circle(.1);
			\node at (3,0) [left] {$a_{2k+2}$};
			\draw (6,0) circle(.1);
			\node at (6,0) [left] {$a_{4k+3}$};
			\node at (7.5,0) [left]{$\cdots$};
			\draw (7.5,0) circle(.1);
			\draw (10.5,0) circle(.1);
			\node at (10.6,0) [left]{$a_{(2k+1)l+1}$};
			
			\draw (0.5,1.5) circle(.1);
			\node at (0.5,1.5) [above] {$a_3$};
			\node at (0.5,1.5) [right]{$\cdots\cdot\cdot$};
			\draw (2.5,1.5) circle(.1);
			\node at (2.5,1.5) [above] {$a_{2k+1}$};
			\draw (3.5,1.5) circle(.1);
			\node at (3.5,1.5) [right]{$\cdots\cdot\cdot$};
			\draw (5.5,1.5) circle(.1);
			\node at (5.5,1.5) [above] {$a_{4k+2}$};	
			\draw ( 8,1.5) circle(.1);
			\node at (8,1.5) [right]{$\cdots\cdot\cdot$};
			\draw ( 10,1.5) circle(.1);
			\node at (10,1.5) [above] {$a_{(2k+1)l}$};
			
			\draw (0.5,-1.5) circle(.1);
			\node at (0.5,-1.5) [below] {$a_2$};
			\node at (0.5,-1.5) [right]{$\cdots\cdot\cdot$};
			\draw (2.5,-1.5) circle(.1);
			\node at (2.5,-1.5) [below] {$a_{2k}$};
			\draw (3.5,-1.5) circle(.1);
			\node at (3.5,-1.5) [right]{$\cdots\cdot\cdot$};
			\draw (5.5,-1.5) circle(.1);
			\node at (5.5,-1.5) [below] {$a_{4k+1}$};	
			\draw ( 8,-1.5) circle(.1);
			\node at (8,-1.5) [right]{$\cdots\cdot\cdot$};
			\draw ( 10,-1.5) circle(.1);
			\node at (10,-1.5) [below] {$a_{(2k+1)l-1}$};
			
			\draw (-3,0) circle(.1);
			\node at (-3,0) [right] {$b_{2k+2}$};
			\draw (-6,0) circle(.1);
			\node at (-6,0) [right] {$b_{4k+3}$};
			\node at (-6.,0) [left]{$\cdots$};
			\draw (-7.5,0) circle(.1);
			
			\draw (-0.5,1.5) circle(.1);
			\node at (-0.5,1.5) [above] {$b_3$};
			\draw (-2.5,1.5) circle(.1);
			\node at (-2.5,1.5) [above] {$b_{2k+1}$};
			\node at (-2.5,1.5) [right]{$\cdots\cdot\cdot$};
			\draw (-3.5,1.5) circle(.1);
			\draw (-5.5,1.5) circle(.1);
			\node at (-5.5,1.5) [above] {$b_{4k+2}$};
			\node at (-5.5,1.5) [right]{$\cdots\cdot\cdot$};	
			\draw (-8,1.5) circle(.1);
			\draw (-10,1.5) circle(.1);
			\node at (-10,1.5) [above] {$b_{(2k+1)l}$};
			\node at (-10,1.5) [right]{$\cdots\cdot\cdot$};
			
			\draw (-0.5,-1.5) circle(.1);
			\node at (-0.5,-1.5) [below] {$b_2$};
			\draw (-2.5,-1.5) circle(.1);
			\node at (-2.5,-1.5) [below] {$b_{2k}$};
			\node at (-2.5    ,-1.5) [right]{$\cdots\cdot\cdot$};
			\draw (-3.5,-1.5) circle(.1);
			\draw (-5.5,-1.5) circle(.1);
			\node at (-5.5,-1.5) [below] {$b_{4k+1}$};	
			\node at (-5.5,-1.5) [right]{$\cdots\cdot\cdot$};
			\draw (-8,-1.5) circle(.1);
			\draw (-10,-1.5) circle(.1);
			\node at (-10,-1.5) [below] {$b_{(2k+1)l-1}$};		
			\node at (-10,-1.5) [right]{$\cdots\cdot\cdot$};

			\draw[-latex,line width=0.5pt](-10.03,1.41)--(-10.3,0.6);
			\draw[line width=0.5pt](-10.2,0.9)--(-10.47,0.09);	
			
			\draw[-latex,line width=0.5pt](-5.53,1.41)--(-5.8,0.6);
			\draw[line width=0.5pt](-5.7,0.9)--(-5.97,0.09);
			
			\draw[-latex,line width=0.5pt](-2.53,1.41)--(-2.8,0.6);
			\draw[line width=0.5pt](-2.7,0.9)--(-2.97,0.09);	
			
			\draw[-latex,line width=0.5pt](0.03,0.09)--(0.3,0.9);	
			\draw[line width=0.5pt](0.47,1.41)--(0.2,0.6);
			
			\draw[-latex,line width=0.5pt](3.03,0.09)--(3.3,0.9);
			\draw[line width=0.5pt](3.47,1.41)--(3.2,0.6);
			
			\draw[-latex,line width=0.5pt](7.53,0.09)--(7.8,0.9);	
			\draw[line width=0.5pt](7.97,1.41)--(7.7,0.6);
			
			\draw[-latex,line width=0.5pt](-7.53,0.09)--(-7.8,0.9);	
			\draw[line width=0.5pt](-7.7,0.6)--(-7.97,1.41);
			
			\draw[-latex,line width=0.5pt](-3.03,0.09)--(-3.3,0.9);	
			\draw[line width=0.5pt](-3.2,0.6)--(-3.47,1.41);
			
			\draw[-latex,line width=0.5pt](2.53,1.41)--(2.8,0.6);	
			\draw[line width=0.5pt](2.97,0.09)--(2.7,0.9);
			
			\draw[-latex,line width=0.5pt](5.53,1.41)--(5.8,0.6);	
			\draw[line width=0.5pt](5.97,0.09)--(5.7,0.9);
			
			\draw[-latex,line width=0.5pt](10.03,1.41)--(10.3,0.6);		
			\draw[line width=0.5pt](10.47,0.09)--(10.2,0.9);
			
			\draw[-latex,line width=0.5pt](-10.03,-1.41)--(-10.3,-0.6);
			\draw[line width=0.5pt](-10.2,-0.9)--(-10.47,-0.09);	
			
			\draw[-latex,line width=0.5pt](-5.53,-1.41)--(-5.8,-0.6);
			\draw[line width=0.5pt](-5.7,-0.9)--(-5.97,-0.09);
			
			\draw[-latex,line width=0.5pt](-2.53,-1.41)--(-2.8,-0.6);
			\draw[line width=0.5pt](-2.7,-0.9)--(-2.97,-0.09);
			
			\draw[line width=0.5pt](0.47,-1.41)--(0.2,-0.6);
			\draw[-latex,line width=0.5pt](0.03,-0.09)--(0.3,-0.9);
			
			\draw[line width=0.5pt](3.47,-1.41)--(3.2,-0.6);
			\draw[-latex,line width=0.5pt](3.03,-0.09)--(3.3,-0.9);
			
			\draw[line width=0.5pt](7.97,-1.41)--(7.7,-0.6);
			\draw[-latex,line width=0.5pt](7.53,-0.09)--(7.8,-0.9);
			
			\draw[-latex,line width=0.5pt](-7.53,-0.09)--(-7.8,-0.9);
			\draw[line width=0.5pt](-7.7,-0.6)--(-7.97,-1.41);	
			
			\draw[-latex,line width=0.5pt](-3.03,-0.09)--(-3.3,-0.9);
			\draw[line width=0.5pt](-3.2,-0.6)--(-3.47,-1.41);
			
			\draw[line width=0.5pt](2.97,-0.09)--(2.7,-0.9);
			\draw[-latex,line width=0.5pt](2.53,-1.41)	--(2.8,-0.6);
			
			\draw[line width=0.5pt](5.97,-0.09)--(5.7,-0.9);
			\draw[-latex,line width=0.5pt](5.53,-1.41)	--(5.8,-0.6);
			
			\draw[line width=0.5pt](10.47,-0.09)--(10.2,-0.9);
			\draw[-latex,line width=0.5pt](10.03,-1.41)--(10.3,-0.6);
			
		\end{tikzpicture}
		\caption{ A broken $k$-diamond of  length $2l$}
		\label{figure-2}
	\end{center}
\end{figure}

\begin{figure}[htbp]
	\begin{center}
		\begin{tikzpicture}[>=Stealth,scale=0.8]	
			\draw (0,0) circle(.1);
			\node at (0,0) [above] {$a_3$};	
			\draw (1.4,0) circle(.1);
			\node at (1.4,0) [above] {$a_5$};	
			\draw (2.8,0) circle(.1);
			\node at (2.8,0) [above] {$a_7$};	
			\node at (2.9,0) [right]{$\cdots\cdots\cdots\cdot$};
			\draw (5.6,0) circle(.1);
			\node at (5.6,0) [above] {$a_{2k-1}$};
			\draw (7,0) circle(.1);
			\node at (7,0) [above] {$a_{2k+1}$};	
			\draw (7.7,-1.4) circle(.1);
			\node at (7.7,-1.4) [right] {$a_{2k+2}$};	
			
			\draw (-0.7,-1.4) circle(.1);
			\node at (-0.7,-1.4) [left] {$a_1$};	
			\draw (0,-2.8) circle(.1);
			\node at (0,-2.8) [below] {$a_2$};
			\draw (1.4,-2.8) circle(.1);
			\node at (1.4,-2.8) [below] {$a_4$};		
			\draw (2.8,-2.8) circle(.1);
			\node at (2.8,-2.8) [below] {$a_6$};
			\node at (2.9,-2.8) [right]{$\cdots\cdots\cdots\cdot$};
			\draw (5.6,-2.8) circle(.1);
			\node at (5.6,-2.8) [below] {$a_{2k-2}$};
			\draw (7,-2.8) circle(.1);
			\node at (7,-2.8) [below] {$a_{2k}$};
			
			\draw[-latex,line width=0.6pt](0.1,0)--(0.8,0);
			\draw[line width=0.6pt] (0.7,0)--(1.3,0);
			\draw[-latex,line width=0.6pt](1.5,0)--(2.2,0);
			\draw[line width=0.6pt] (2.1,0)--(2.7,0);
			\draw[-latex,line width=0.6pt](5.7,0)--(6.4,0);
			\draw[line width=0.6pt] (6.3,0)--(6.9,0);
			
			\draw[-latex,line width=0.6pt](7.045,-0.09)--(7.4,-0.8);
			\draw[line width=0.6pt] (7.3,-0.6)--(7.655,-1.31);
			\draw[-latex,line width=0.6pt](-0.655,-1.31)--(-0.3,-0.6);
			\draw[line width=0.6pt] (-0.4,-0.8)--(-0.045,-0.09);
			\draw[-latex,line width=0.6pt](-0.655,-1.49)--(-0.3,-2.2);
			\draw[line width=0.6pt] (-0.4,-2)--(-0.045,-2.71);
			
			\draw[-latex,line width=0.6pt](0.1,-2.8)--(0.8,-2.8);
			\draw[line width=0.6pt] (0.7,-2.8)--(1.3,-2.8);
			\draw[-latex,line width=0.6pt](1.5,-2.8)--(2.2,-2.8);
			\draw[line width=0.6pt] (2.1,-2.8)--(2.7,-2.8);
			\draw[-latex,line width=0.6pt](5.7,-2.8)--(6.4,-2.8);
			\draw[line width=0.6pt] (6.3,-2.8)--(6.9,-2.8);
			\draw[-latex,line width=0.6pt](7.045,-2.71)--(7.4,-2);
			\draw[line width=0.6pt] (7.3,-2.2)--(7.655,-1.49);
			
			\draw[-latex,line width=0.6pt](0.045,-0.09)--(0.4,-0.8);
			\draw[line width=0.6pt] (0.3,-0.6)--(1.355,-2.71);
			\draw[-latex,line width=0.6pt](1.445,-0.09)--(1.8,-0.8);
			\draw[line width=0.6pt] (1.7,-0.6)--(2.755,-2.71);
			\draw[-latex,line width=0.6pt](5.645,-0.09)--(6,-0.8);
			\draw[line width=0.6pt] (5.9,-0.6)--(6.955,-2.71);
			
			\draw[-latex,line width=0.6pt](0.045,-2.71)--(0.4,-2);
			\draw[line width=0.6pt] (0.3,-2.2)--(1.355,-0.09);
			
			\draw[-latex,line width=0.6pt](1.445,-2.71)--(1.8,-2);
			\draw[line width=0.6pt] (1.7,-2.2)--(2.755,-0.09);
			
			\draw[-latex,line width=0.6pt](5.645,-2.71)--(6,-2);
			\draw[line width=0.6pt] (5.9,-2.2)--(6.955,-0.09);
			
		\end{tikzpicture}
		
		\caption{ A $k$-elongated partition diamond of  length $1$}
		\label{figure-1}
	\end{center}
\end{figure}
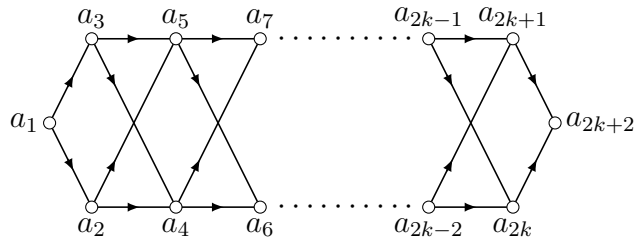

The main objective of this paper is to establish Tur\'an inequalities for the number of broken $k$-diamond partitions. The notion of broken $k$-diamond partitions was introduced by Andrews and Paule \cite{Andrews-Paule-2007}. A broken $k$-diamond partition $\pi=(b_2,\,\ldots, b_{2k+2},\ldots, b_{(2k+1)l+1},\,a_1,\,\ldots,a_{2k+2}, \ldots, a_{(2k+1)l+1})$ is a plane partition satisfying the relations illustrated in Figure \ref{figure-2},  where $a_i,\,b_i$ are nonnegative integers and $a_i\rightarrow a_j$ is interpreted as $a_i\geq a_j$. More precisely, each building block in Figure \ref{figure-2}  has the same order structure as shown in Figure \ref{figure-1}. We call such block a $k$-elongated partition diamond of length $1$. It should be noted that the broken block $(b_2,b_3,\ldots, b_{2k+2})$ is also a $k$-elongated partition diamond of length $1$ from which a source $b_1$ is deleted.

Let $\Delta_k(n)$ denote the number of broken $k$-diamond partitions of $n$. Andrews and Paule \cite{Andrews-Paule-2007} established the following generating function for $\Delta_k(n)$:
\begin{align*}
\sum_{n=0}^{\infty} \Delta_{k}(n) q^{n} &=\prod_{n=1}^{\infty} \frac{\left(1-q^{2 n}\right)\left(1-q^{(2 k+1) n}\right)}{\left(1-q^{n}\right)^{3}\left(1-q^{(4 k+2) n}\right)}.
\end{align*}
It's known that $\Delta_{k}(n)$ possesses many beautiful arithmetic properties. Many Ramanujan-like congruences satisfied by $\Delta_1(n)$ and $\Delta_2(n)$  were proved by Andrews and Paule \cite{Andrews-Paule-2007}  and other authors, see, for example, Chan \cite{Chan-2008}, Chen, Fan and Yu \cite{Chen-Fan-Yu-2014}, Hirschhorn \cite{Hirschhorn-2018},  Paule and Radu \cite{Paule-Radu-2010} and so on. It should be noted that  $\Delta_k(n)$ are the coefficients of a  modular function  with respect to $\Gamma_0(4k+2)$.

The Tur\'an inequalities (or sometimes called Newton inequalities) arise in the study of
real entire functions in Laguerre-P\'olya class, which are  closely related to the study of the
Riemann hypothesis \cite{Dimitrov-1998, Szego-1948}. A sequence $\{\alpha_n\}_{n\geq 0}$ of real numbers is log-concave if it satisfies the (second order) Tur\'an inequalities $\alpha^2_n\geq \alpha_{n-1}\alpha_{n+1}$ for $n\geq 1$. We call the sequence $\{\alpha_n\}_{n\geq 0}$ satisfies the third order Tur\'an inequalities if for $n \geq 1$,
\[4(\alpha^2_n-\alpha_{n-1}\alpha_{n+1})(\alpha^2_{n+1}-\alpha_n\alpha_{n+2})\geq (\alpha_{n}\alpha_{n+1}-\alpha_{n-1}\alpha_{n+2})^2.\]

As stated by Chen, Jia and Wang \cite{Chen-Jia-Wang-2019} and Griffin, Ono, Rolen and Zagier \cite{Griffin-Ono-Rolen-Zagier-2019}, the higher order Tur\'an inequalities are conveniently formulated in terms of the Jensen polynomials. The  Jensen polynomials $J_{\alpha}^{d,n}(X)$ of degree $d$ and shift $n$ associated to the sequence $\{\alpha_n\}_{n\geq0}$ are defined by
\begin{equation}
	J_{\alpha}^{d,n}(X)=\sum_{i=0}^{d} \binom{d}{i} \alpha_{n+i} X^i.
\end{equation}
When $d=2$ and shift $n-1$, the Jensen polynomial $J_{\alpha}^{2,n-1}(X)$ reduces to
\[J_{\alpha}^{2,n-1}(X)=\alpha_{n-1}+2\alpha_{n} X+\alpha_{n+1} X^2.
\]
It is clear that $\{\alpha_n\}_{n\geq0}$ is log-concave at $n$ if and only if  $J_{\alpha}^{2,n-1}(X)$ has only real roots. In general, we say that the sequence $\{\alpha_n\}_{n\geq0}$ satisfies the order $d$ Tur\'an inequality at $n$ if and only if $J_{\alpha}^{d,n-1}(X)$ is hyperbolic. Recall that a polynomial is hyperbolic if all of its roots are real.

There are several recent work on the Tur\'an inequalities for the partition functions. Nicolas \cite{Nicolas-1978} and DeSalvo and Pak \cite{DeSalvo-Pak-2015} proved that the partition function $p(n)$ is log-concave
for $n \geq 26$, where $p(n)$ is the number of partitions of $n$. Chen \cite{Chen-2017}  conjectured that $p(n)$ satisfies the third order Tur\'an inequalities
for $n \geq 95$, which was proved by Chen, Jia and Wang \cite{Chen-Jia-Wang-2019}.  Chen, Jia and Wang \cite{Chen-Jia-Wang-2019} further conjectured that for $d\geq 4$, there exists a
positive integer $N_p(d)$ such that $p(n)$ satisfies the order $d$  Tur\'an inequalities  for $n\geq N_p(d)$, that is,   the Jensen polynomial $J_{p}^{d,n-1}(X)$ associated to $p(n)$  is hyperbolic for $n\geq N_p(d)$. Griffin, Ono, Rolen and Zagier \cite{Griffin-Ono-Rolen-Zagier-2019} showed that their conjecture is true for   sufficiently large $n$. In fact,  they obtained the following general result:

\begin{thm}[Proof of Theorem 7 of \cite{Griffin-Ono-Rolen-Zagier-2019} ] \label{Ono} Let $\{a_f(n)\}_{n\geq 0}$ be a sequence of positive real numbers. If
\[a_f(n)=A_f n^{\frac{k-1}{2}}I_{k-1}(4\pi \sqrt{mn})+O(n^Ce^{2\pi\sqrt{mn}})\]
as $n\rightarrow \infty$ for some non-zero constants $A_f, m$ and $C$, where $I_{\nu}(s)$    is the $\nu$-th modified Bessel function of the first kind. Then for $d\geq 1$, the Jensen polynomial $J_{a_f}^{d,n}(X)$ associated to $a_f(n)$ is hyperbolic for all sufficiently large $n$.
\end{thm}

Since then,  Tur\'an inequalities for  other partition functions have been extensively investigated. For example, Engel \cite{Engel-2017} showed that the overpartition function $\overline{p}(n)$ is log-concave
for $n \geq 2$ and Liu and Zhang \cite{Liu-Zhang-2021} showed that the overpartition function $\overline{p}(n)$  satisfies the third order Tur\'an inequalities
for $n \geq 16$. Recently, Bringmann, Kane, Rolen and Tripp \cite{Bringmann-Kane-Rolen-Tripp-2021} showed that $k$-colored partition function $p_k(n)$ is log-concave for $n\geq 6$.    Craig and Pun \cite{Craig-Pun-2021} showed that  the number of the $k$-regular partitions of $n$ satisfies   the order $d\geq 1$  Tur\'an inequalities  for sufficiently large $n$ and Ono, Pujahari and Rolen \cite{Ono-Pujahari-Rolen-2019}  showed that the number of the number of MacMahon's plane partitions of $n$ satisfies   the order $d\geq 1$  Tur\'an inequalities  for sufficiently large $n$. It should be noted that  Craig and Pun's result can be viewed as a direct consequence of  Theorem \ref{Ono}.

In this paper, we intend to explore Tur\'an inequalities for  broken $k$-diamond partitions.   Appealing to Sussman's Rademacher-type formula for $\eta$-quotients \cite{Sussman-2017} , we obtain the following asymptotic formula for $\Delta_k(n)$, where $k=1$ or $2$.

\begin{thm}\label{asympform} For $k=1$ or $2$, as $n\rightarrow \infty$,
\begin{equation}\label{eq-del-asy}
	\Delta_k(n)=\frac{(5k+2) \pi^3}{18 (2k+1) {x^2_k(n)}} I_2\left(\sqrt{\frac{5k+2}{2k+1} } {x_k(n)}\right) +O\left(x^{-7/2}_k(n)\exp{\left(\frac{\sqrt{ {5k+2}}}{2\sqrt{2k+1}}  {x_k(n)}\right)}\right),
\end{equation}
  where $I_{2}(s)$ is the second  modified Bessel function of the first kind, and
\begin{equation}\label{defi-xk}
x_k(n)=\frac{\pi\sqrt{24n-(2k+2)}}{6}.
\end{equation}

\end{thm}

Combining Theorem \ref{Ono} and Theorem \ref{asympform},  we derive that for $k=1$ or $2$ and $d\geq 1$, $\Delta_k(n)$ satisfies the order $d$ Tur\'an inequalities for  sufficiently large $n$. To wit,

\begin{coro}\label{the-gen}
Let $\Delta_k=\{\Delta_k(n)\}_{n\geq 0}$, for $k=1,2$ and $d\geq 1$, the Jensen polynomial $J_{\Delta_k}^{d,n}(X)$ associated to $\Delta_k$ is hyperbolic for all sufficiently large $n$.
\end{coro}

It is worth mentioning that   there exists a minimal natural number $N_{\Delta_k}(d)$ such that  $J_{\Delta_k}^{d,n}(X)$ is hyperbolic for $n\geq N_{\Delta_k}(d)$. Table \ref{conj-value} provides conjectural  values for $N_{\Delta_{k}}(d)$  for   $ 2\leq d\leq 13$.

\begin{table}[h]
\centering
	\begin{tabular}{l|l|l|l|l|l|l|l|l|l|l|l|l}\hline
	$d$ & $2$ &$3$ & 4&5 &6&7&8&9&10&11&12&13\\\hline
	$N_{\Delta_1}(d)$ & $0$  & 4 & 17 & 41&72&116&171&238&320&415&525&650\\\hline
	$N_{\Delta_2}(d)$ &$0$ & 4 & 17 & 34&62&99&147&200&272&355&445&552\\
		\hline
	\end{tabular}
	\caption{The conjectural  values of $N_{\Delta_k}(d)$  for $1\leq k\leq 2$ and  $2\leq d\leq 13$.} \label{conj-value}
\end{table}

We further  prove that $N_{\Delta_k}(2)=0$ where $k=1$ or $2$. More precisely, we show that

\begin{thm}\label{thm-Delta-1}
	For $k=1$ or $2$, the broken $k$-diamond partition function $\Delta_k(n)$ is log-concave for $n\geq 1$, that is,
	\begin{equation}\label{thm-Delta-1-eq}
		\Delta_k^2(n)\geq \Delta_k(n-1)\Delta_k(n+1).
	\end{equation}
\end{thm}

As noted in the paper by Asai, Kubo and Kuo \cite{Asai-Kubo-Kuo-00} and Sagan \cite{Sagan-1988}, we see that Theorem \ref{thm-Delta-2} is equivalent to the following multiplicative properties of $\Delta_k(n)$.

\begin{coro}\label{thm-Delta-2}
	For $k=1$ or $2$ and   $a,\,b\ge 1$,
	\begin{equation*}
	\Delta_k(a)\Delta_k(b)\ge\Delta_k(a+b).
	\end{equation*}
\end{coro}
It should be noted that  the multiplicative properties of the ordinary partition function $p(n)$ were initially obtained by Bessenrodt and Ono \cite{Bessenrodt-Ono-2016}. Subsequently, the multiplicative properties of other partition functions have been established, for example,
Beckwith and Bessenrodt  \cite{Beckwith-Bessenrodt-2016} established the multiplicative properties of  $k$-regular partition function  and   Bringmann, Kane, Rolen and Tripp \cite{Bringmann-Kane-Rolen-Tripp-2021}  acquired the   multiplicative properties of $k$-colored partition function, which resolved a conjecture of  Chern, Fu and Tang \cite{Chern-Fu-Tang-2018}.

This article is organized as follows. In Section 2, we prove Theorem \ref{asympform} with the aid of  Sussman's Rademacher-type formula for $\eta$-quotients. Section 3 is devoted to the proof of Theorem \ref{thm-Delta-1}. To this end, we  derive an upper bound and a lower bound  for $\Delta_k(n)$ with the aid of Theorem \ref{asympform} and establish  an inequality on the second Bessel function.   In Section 4,  we pose some questions and remarks for future work.

\section{Proof of Theorem \ref{asympform}}

To prove Theorem \ref{asympform}, we first derive Rademacher-type formulas for $\Delta_k(n)$ ($k=1$ or $2$) with the aid of Sussman's Rademacher-type formula for $\eta$-quotients  \cite{Sussman-2017}.  Define
	\begin{eqnarray}
	G(q):=\prod_{r=1}^R(q^{m_r};q^{m_r})_{\infty}^{\delta_r},
	\label{eq-1-1}
	\end{eqnarray}
where $\mathbf{m}=(m_1,\ldots,m_R)$ is a sequence of  $R$ distinct positive integers and $\mathbf{\delta}=(\delta_1,\ldots,\delta_R)$ is a sequence of $R$ non-zero integers.
Here and throughout this paper,  we have adopted the  standard notation on  $q$-series \cite{Andrews-1998}:
\[(a;q)_n=\prod_{j=0}^{n-1}(1-aq^j) \quad \text{and} \quad (a;q)_\infty=\prod_{j=0}^\infty (1-aq^j).\]
In order to present Sussman's result, we need a few preliminary definitions. Let
\[n_0=-\sum_{r=1}^Rm_r\delta_r.\]
The function $q^{\frac{n_0}{24}}G(q)$ is holomorphic in the open unit disk $D$, and so we may write
 \[G(q)=q^{-\frac{n_0}{24}}\sum_{n\ge 0}g(n)q^n\]
 for some coefficients $g(n)$. Sussman \cite{Sussman-2017} obtained a Rademacher-type formula for $g(n)$ with $\frac{1}{2}\sum_{r=1}^R\delta_r<0$, which is a special case of   Bringmann and Ono \cite{Bringmann-Ono-2012}.  Let
\[c_1=-\frac{1}{2}\sum_{r=1}^R\delta_r, \quad c_2(j)=\prod_{r=1}^R\left(\frac{\gcd(m_r,j)}{m_j}\right)^{\frac{\delta_r}{2}},\quad
c_3(j)=-\sum_{r=1}^R\frac{\delta_r\gcd^2(m_r,j)}{m_r},\]
\begin{equation}\label{defi-A}
\hat{A}_{j}(n)=\sum_{0\leq h< j\atop {\rm gcd}(h,j)=1}\exp\left(-\frac{2\pi h i }{j}-\pi i\sum_{r=1}^R\delta_rs\left(\frac{m_rh}{\gcd(m_r,j)},\frac{j}{\gcd(m_r,j )}\right)\right),
\end{equation}
where  $s(h,j)$ is the  Dedekind sum:
\[s(h,j)=\sum_{r=1}^{j-1}\left(\frac{r}{j}-\left[\frac{r}{j}\right]-\frac{1}{2}\right)\left(\frac{h r}{j}-\left[\frac{h r}{j}\right]-\frac{1}{2}\right).\]

\begin{thm}(Sussman)\label{the-sussman}
If  $c_1>0$ and the periodic function $\beta(j) \colon \mathbb{N} \rightarrow \mathbb{R}$ given by
\begin{equation}
	\beta(j)=\min_{1\leq r \leq R}{\frac{{\rm gcd}^2(m_r,j)}{m_r}}-\frac{c_3(j)}{24}
\end{equation}
is non-negative, then for  $n>\frac{n_0}{24}$, we have
		\begin{align}
		g(n)&=2\pi \left(\frac{1}{{24n-n_0}}\right)^{\frac{c_1+1}{2}} \sum_{j\geq 1\atop c_3(j)\geq 0}  c_2(j) c_3(j )^{\frac{c_1+1}{2}}j^{-1}\hat{A}_j(n){I_{c_1+1}\left(\frac{\pi \sqrt{c_3(j)(24n-n_0)}}{6j}\right)},
		\end{align}	
		where
 $I_{\nu}(s)$ is the $\nu$-th modified Bessel function of the first kind.
	\end{thm}
Here and throughout this paper,  we  adopt  the following infinite series definition of   the  modified Bessel function of the first kind,
\[I_{\nu}(s):=\sum_{r=0}^{\infty} \frac{1}{r!\Gamma(r+\nu+1)}\left(\frac{s}{2}\right)^{2 r+\nu}.\]

By invocation of Theorem \ref{the-sussman}, we attain the following Rademacher-type formula for $\Delta_k(n)$ ($k=1$ or $2$).

\begin{thm}\label{eq-Del-k-thm}For $k=1$ or $2$ and $n\geq 1$,
\begin{equation}\label{eq-Del-k}
	\Delta_k(n)=\frac{   \pi^3}{18 x^2_k(n)}\sum_{j\geq 1 } \alpha_k(j) j^{-1} \hat{A}_j(n){I_{2}\left(\frac{\sqrt{\alpha_k(j)} x_k(n)}{j}\right)},
\end{equation}
where $x_k(n)$ is defined in \eqref{defi-xk}, $\hat{A}_j(n)$ is defined in \eqref{defi-A},   $I_2$ is the second modified Bessel function of the first kind, and
\begin{equation}\label{defi-alphak}
\alpha_k(j):=\left\{
\begin{aligned}
	&1+\frac{{\rm gcd}^2(2k+1,j)}{2k+1}, \quad  \text{$j$ is even},\\[6pt]
	&\frac{5}{2}-\frac{{\rm gcd}^2(2k+1,j)}{4k+2}, \quad \text{$j$ is odd}.
\end{aligned}
\right.
\end{equation}

\end{thm}
\proof Recall that	
\begin{align*}
\sum_{n=0}^{\infty} \Delta_{k}(n) q^{n} =  \frac{(q^2;q^2)_\infty(q^{2k+1};q^{2k+1})_\infty}
{(q;q)_\infty^{3}(q^{4k+2};q^{4k+2})_\infty}.
\end{align*}
We have
\[\mathbf{m}=(1,2,2k+1,4k+2) \quad \text{and} \quad \mathbf{\delta}=(-3,1,1,-1),\]
so that  $n_0=2k+2$, $c_1=1,$ and for $j\geq 1$,
\begin{align*}
c_2(j)&=\left(\frac{\gcd(2,j)}{2}\right)^{\frac{1}{2}}
\left(\frac{\gcd(2k+1,j)}{2k+1}\right)^{\frac{1}{2}}
\left(\frac{\gcd(4k+2,j)}{4k+2}\right)^{-\frac{1}{2}}\\[5pt]
&=\left(\frac{\gcd(2,j)\gcd(2k+1,j)}{\gcd(4k+2,j)}\right)^{\frac{1}{2}}=1,
\end{align*}
\begin{align}
c_3(j)&=3-\frac{\gcd^2(2,j)}{2}-\frac{\gcd^2(2k+1,j)}{2k+1}
+\frac{\gcd^2(4k+2,j)}{4k+2}\nonumber\\[5pt]
&=\left\{
\begin{aligned}
	&1+\frac{\gcd^2(2k+1,j)}{2k+1}, \quad  \text{$j$ is even},\\[6pt]
	&\frac{5}{2}-\frac{ \gcd^2(2k+1,j)}{4k+2}, \quad \text{$j$ is odd},
\end{aligned}
\right. \label{c3j}
\end{align}
and
\begin{align}
\min_{1\leq r\leq 4}\left(\frac{\gcd^2(m_r,j)}{m_r}\right)&=\min\left\{1,
\frac{\gcd^2(2,j)}{2},\frac{\gcd^2(2k+1,j)}{2k+1},
\frac{\gcd^2(2k+1,j)\gcd^2(2,j)}{4k+2}\right\}\nonumber\\[9pt]
&=\left\{
\begin{aligned}
	&\min\left\{1,\frac{\gcd^2(2k+1,j)}{2k+1}\right\}, \quad  \text{$j$ is even},\\[9pt]
	&\min\left\{\frac{1}{2},\frac{\gcd^2(2k+1,j)}{4k+2} \right\}, \quad \text{$j$ is odd}.   \label{condit}
\end{aligned}
\right.
\end{align}
Combining \eqref{c3j} and  \eqref{condit}, we get
\begin{align*}
\beta_k(j)=
\left\{
\begin{aligned}
	&\min\left\{1,\frac{\gcd^2(2k+1,j)}{2k+1}\right\}-\frac{1}{24}\left(1+\frac{\gcd^2(2k+1,j)}{2k+1}\right), \quad  \text{$j$ is even},\\[9pt]
	&\min\left\{\frac{1}{2},\frac{\gcd^2(2k+1,j)}{4k+2} \right\}-\frac{1}{24}\left(\frac{5}{2}-\frac{ \gcd^2(2k+1,j)}{4k+2}\right), \quad \text{$j$ is odd}.
\end{aligned}
\right.
\end{align*}
Set $\alpha_k(j)=c_3(j)$. From the definitions of $\alpha_k(j)$ and $\beta_k(j)$, we see that  $\alpha_k(j)$ and $\beta_k(j)$ are periodic functions with period $4k+2$. The following table gives the values of $\alpha_k(j)$ and $\beta_k(j)$ for $1\leq  j\leq 4k+2$.
\begin{table}[h]
	\centering
	\begin{tabular}{l|l|l|l|l|l|l|l|l|l|l|l}\hline
	
		$j$&$1$&$2$&$3$&$4$&$5$&$6$&$7$&$8$&$9$&$10$&$\cdots$\\ \hline
		
		$\alpha_1(j)$& $\frac{7}{3}$&$\frac{4}{3}$&$1$&$\frac{4}{3}$ &$\frac{7}{3}$&$4$& $\frac{7}{3}$&$\frac{4}{3}$&$1$&$\frac{4}{3}$&$\cdots$\\\hline
		
		$\alpha_2(j)$&$\frac{12}{5}$&$\frac{6}{5}$&$\frac{12}{5}$&$\frac{6}{5}$&$0$
&$\frac{6}{5}$&$\frac{12}{5}$&$\frac{6}{5}$&$\frac{12}{5}$&$6$&$\cdots$\\
		\hline
		$\beta_1(j)$ &$\frac {5} {72}$ &$\frac {5} {18} $&$\frac {11} {24} $&$\frac {5} {18} $&$\frac {5} {72}$ &$\frac {5} {6}$ &$\frac {5} {72}$ &$\frac {5} {18}$ &$\frac{11} {24} $&$\frac {5} {18} $& $\cdots$\\\hline
$\beta_2(j)$&$0$ &$\frac {3} {20}$ &$ 0$ &$\frac {3} {20}$ &$\frac {1} {2}$ &$\frac {3} {20}$ & $0$ &$\frac {3} {20}$ &$ 0 $&$\frac {3} {4}$&$\cdots$ \\\hline
	\end{tabular}
	\caption{The values of $\alpha_k(j)$ and $\beta_k(j)$ for $k=1$ or $k=2$. }\label{table-sec2}
\end{table}

From Table \ref{table-sec2}, we find that  $\alpha_k(j)\geq0$ when $k=1$ or $2$ and $\Delta_k(n)$ satisfies the condition  in Theorem \ref{the-sussman} when $k=1$ or $2$.   Hence we derive \eqref{eq-Del-k} by substituting relevant values  into Theorem \ref{the-sussman}. This completes the proof. \qed

It should be noted that Theorem \ref{the-sussman} could not be applied to derive the explicit formula for $\Delta_k(n)$ when $k\geq 3$.  Setting $j=1$, we find that
\[\min_{1\leq r\leq 4}\left(\frac{\gcd^2(m_r,j)}{m_r}\right)=\min\left\{\frac{1}{2},\frac{1}{4k+2} \right\}=\frac{1}{4k+2},\]
and
\[c_3(1)=\frac{5}{2}-\frac{1}{4k+2}.\]
But when $k\geq 3$,
\[\frac{1}{4k+2}<\frac{1}{24}\left(\frac{5}{2}-\frac{1}{4k+2}\right),\]
which implies that $\Delta_k(n)$ does not satisfy  the condition in Theorem \ref{the-sussman} when $k\geq 3$.  Here and in the sequel, we assume that $k=1$ or $2$.

We are now in a position to prove Theorem \ref{asympform}  by means of  Theorem \ref{eq-Del-k-thm}.

\noindent{\it Proof of Theorem \ref{asympform}.} Define
\begin{equation}\label{defi-Mk}
M_k(n):=\frac{\alpha_k(1) \pi^3}{18  {x^2_k(n)}}  I_2\left(\sqrt{\alpha_k(1) } {x_k(n)}\right).
\end{equation}

Observing that $\hat{A}_1(n)=1$ and $\alpha_k(1)=\frac{5k+2}{2k+1}$, we deduce from Theorem \ref{eq-Del-k-thm} that
\begin{equation}\label{fn1}
 \Delta_k(n)= M_k(n)+ R_k(n),
\end{equation}
where
\[R_k(n)=\frac{ \pi^3}{18 x^2_k(n)}\sum_{j\geq 2 } \alpha_k(j) j^{-1} \hat{A}_j(n){I_{2}\left(\frac{\sqrt{\alpha_k(j)} x_k(n)}{j}\right)}.\]
We next establish the upper bound for $|R_k(n)|$.

By the definition of  $\hat{A}_j(n)$, we derive that  for any $n\geq 0$ and $j\geq 1$,
 \begin{equation}\label{eq-A-ineq}
 |\hat{A}_j(n)|\leq j,
 \end{equation}
since $|e^{2\pi r i}|=1$ for any $r\in \mathbb{R}$.

In light of the fact that $\alpha_k(j)$ is a periodic function with period $4k+2$, we see from Table \ref{table-sec2} that
 \begin{equation}\label{eq-max-k}
 	 \max_{ j\geq 2\atop k=1,2} \alpha_k(j)=\max_{2\leq j\leq 4k+2\atop k=1,2} \alpha_k(j)=6.
 \end{equation}
Combing \eqref{eq-A-ineq} and \eqref{eq-max-k}, we are led to
\begin{align}
\label{eq-R-k-1} |R_k(n)|&\leq \frac{ \pi^3}{18 x^2_k(n)}\sum_{j\geq 2 } |\alpha_k(j)|j^{-1} |\hat{A}_j(n)|{I_{2}\left(\frac{\sqrt{\alpha_k(j)} x_k(n)}{j}\right)}\\[3pt]
	&\leq \frac{\pi^3}{3x^2_k(n)}\sum_{j\geq 2} I_{2}\left(\frac{\sqrt{\alpha_k(j)} x_k(n)}{j}\right) 	\nonumber\\[3pt]
\nonumber	&=\frac{\pi^3}{3x^2_k(n)} \sum_{j\geq 2 \atop j\nmid 4k+2}I_{2}\left(\frac{\sqrt{\alpha_k(j)} x_k(n)}{j}\right)+\frac{\pi^3}{3x^2_k(n)} \sum_{j\geq 2 \atop j\mid 4k+2}I_{2}\left(\frac{\sqrt{\alpha_k(4k+2)} x_k(n)}{j}\right) 	\nonumber\\[3pt]
	&=\frac{\pi^3}{3x^2_k(n)} \sum_{j\geq 2 \atop j\nmid 4k+2}I_{2}\left(\frac{\sqrt{\alpha_k(j)} x_k(n)}{j}\right)+\frac{\pi^3}{3x^2_k(n)} \sum_{j\geq 1}I_{2}\left(\frac{\sqrt{\alpha_k(4k+2)} x_k(n)}{j(4k+2)}\right). \nonumber
\end{align}
It is evident that
\begin{equation*}
	\max_{j\geq 2 \atop j\nmid 4k+2} \alpha_k(j)=\max_{2\leq j<4k+ 2}\alpha_k(j)\leq \alpha_k(1),
\end{equation*}
 and
\begin{equation*}	\frac{\sqrt{\alpha_k(4k+2)}}{4k+2}\leq\frac{\sqrt{\alpha_k(1)}}{2}.
\end{equation*}
Hence \eqref{eq-R-k-1} can be further bounded by
\begin{equation*}
	|R_k(n)|\leq  \frac{\pi^3}{3x^2_k(n)} \sum_{j\geq 2}I_{2}\left(\frac{\sqrt{\alpha_k(1)}x_k(n)}{j}\right)+\frac{\pi^3}{3x^2_k(n)} \sum_{j\geq 1}I_{2}\left(\frac{\sqrt{\alpha_k(1)} x_k(n)}{2j}\right),
\end{equation*}
since  $I_2(s)$ is increasing on $(0, \infty)$.

Note that
\begin{align*}
	\sum _{j\geq N} I_2\left(\frac{s}{j}\right) & =\sum _{j\geq N} \sum _{m\geq 0} \frac{1}{m! (m+2)!}\left(\frac{s}{2 j}\right)^{2 m+2}\\
	&\leq \int _N^{\infty }\sum _{m\geq 0} \frac{1}{m! (m+2)!} \left(\frac{s}{2 t}\right)^{2 m+2} {\rm d}t\\
	&=\sum _{m\geq 0} \frac{1}{m! (m+2)!} \int_N^{\infty } \left(\frac{s}{2 t}\right)^{2 m+2} {\rm d}t\\
	%&=\sum_{m\geq 0} \frac{t^{2m+2}}{2^{m+2}N^{2m+1}(2m+1)m!(m+2)!}\\
	&=N\sum _{m\geq 0 } \frac{ 1}{(2 m+1) m! (m+2)!} \left(\frac{s}{2N}\right)^{2 m+2}\\
	&\leq N \sum _{m\geq  0} \frac{1}{(m+1)! (m+2)!}\left(\frac{s}{2 N}\right)^{2 m+2}\\
	&=N\sum _{m\geq  1} \frac{1}{m! (m+1)!}\left(\frac{s}{2 N}\right)^{2 m}\\
	&\leq  \frac{2 N^2}{s}  \sum _{m\geq 0} \frac{1}{m!(m+1)!} \left(\frac{s}{2 N}\right)^{2 m+1}\\
	&=  \frac{2N^2}{s}I_1\left(\frac{s}{N}\right).
\end{align*}
Thus we obtain that
\begin{align}
	\nonumber|R_k(n)|&\leq  \frac{\pi^3}{3x^2_k(n)} \sum_{j\geq 2}I_{2}\left(\frac{\sqrt{\alpha_k(1)}x_k(n)}{j}\right)
+\frac{\pi^3}{3x^2_k(n)} \sum_{j\geq 1}I_{2}\left(\frac{\sqrt{\alpha_k(1)} x_k(n)}{2j}\right)\nonumber\\[6pt]
	\nonumber &\leq \frac{8\pi^3}{3\sqrt{\alpha_k(1)}x^3_k(n)} I_{1}\left(\frac{1 }{2}\sqrt{\alpha_k(1)}x_k(n) \right)+\frac{4\pi^3}{3\sqrt{\alpha_k(1)}x^3_k(n)} I_{1}\left(\frac{1 }{2}\sqrt{\alpha_k(1)}x_k(n)\right)
	\nonumber \\[6pt]
	\nonumber &\leq \frac{4\pi^3}{\sqrt{\alpha_k(1)} x_k^3(n)} I_{1}\left(\frac{1 }{2}\sqrt{\alpha_k(1)}x_k(n)\right).
\end{align}
Using Lemma 2.2 (1) of Bringmann, Kane, Rolen and Tripp \cite{Bringmann-Kane-Rolen-Tripp-2021}, we find that for $s\geq 1$,
\begin{equation*}
I_1(s) \leq \sqrt{\frac{2}{\pi s}} e^s.
\end{equation*}
Hence
\begin{equation}\label{bound-Rk}
	|R_k(n)|\leq \frac{8 \pi^{\frac{5}{2}}}{ {\alpha^{\frac{3}{4}}_k(1)} x^{\frac{7}{2}}_k(n)}  \exp\left(\frac{1}{2}\sqrt{\alpha_k(1)}x_k(n)\right) .
\end{equation}
Consequently,  we derive \eqref{asympform} from \eqref{fn1} and \eqref{bound-Rk} upon noting that
\[x_k(n)=\frac{\pi\sqrt{24n-(2k+2)}}{6}\rightarrow \infty \quad \text{as} \quad n\rightarrow \infty.\]  This completes the proof of Theorem \ref{asympform}. \qed

\section{Proof of  Theorem \ref{thm-Delta-1} }

To prove Theorem \ref{thm-Delta-1},  we    establish an upper bound and a lower bound  for $\Delta_k(n)$  in light of Theorem \ref{eq-Del-k-thm} and   an inequality on the second Bessel function.

\begin{thm}\label{lem-Delta-1} Let $x_k(n)$ be defined as in \eqref{defi-xk}, let $\alpha_k(n)$ be defined as in \eqref{defi-alphak} and let $M_k(n)$ be defined as in \eqref{defi-Mk}. For $k=1$ or $2$ and $x_k(n) \geq 152 $, we have
	\begin{equation}\label{eq-main}
		M_k(n)\left(1-\frac{1}{x^6_k(n)}\right)\leq \Delta_k(n)\leq M_k(n)\left(1+\frac{1}{x^6_k(n)}\right).
	\end{equation}
	
\end{thm}

{\noindent \it Proof.} Define
\begin{equation}\label{defi-Gk}
	G_k(n):=\frac{\frac{8 \pi^{\frac{5}{2}}}{ {\alpha^{\frac{3}{4}}_k(1)} x^{\frac{7}{2}}_k(n)}  \exp\left(\frac{1}{2}\sqrt{\alpha_k(1)}x_k(n)\right)}
	{\frac{ {\alpha_k(1)}\pi^3}{18  x^2_k(n)}  I_2\left(\sqrt{\alpha_k(1) } x_k(n)\right)}=\frac{144}{\alpha_k^{\frac{7}{4}}(1)\sqrt{\pi}x_k^{\frac{3}{2}}(n)} \cdot
	\frac{ \exp\left(\frac{1}{2}\sqrt{\alpha_k(1)}x_k(n)\right)}
	{ I_2\left(\sqrt{\alpha_k(1) } x_k(n)\right)}.
\end{equation}
Based on  \eqref{bound-Rk}, we derive that
\begin{equation*}
	M_k(n)(1-G_k(n))\leq \Delta_k(n)\leq M_k(n)(1+G_k(n)).
\end{equation*}
To show \eqref{eq-main}, it is enough to prove that for $x_k(n)\geq 152$,
\begin{equation}
	G_k(n)\leq \frac{1}{x^6_k(n)}.
\end{equation}
Invoking   Lemma 2.2 (4) of Bringmann, Kane, Rolen and Tripp \cite{Bringmann-Kane-Rolen-Tripp-2021}, we know that for $s\geq 231$,
\begin{equation}\label{eq-BesselI1}
	\left|I_2(s) e^{-s} \sqrt{2\pi s}-1+\frac{15}{8 s}-\frac{105}{128 s^2}-\frac{315}{1024 s^3}  \right|\leq \frac{3968}{3 s^4}.
\end{equation}
Hence, for $s\geq 231$,
\begin{equation*}
	I_2(s)\geq \frac{e^s}{\sqrt{2\pi s}}\left(1-\frac{15}{8s}+\frac{105}{128s^2}+\frac{315}{1024s^3}
-\frac{3968}{3s^4}\right).
\end{equation*}
Note that
\begin{equation*}
	\frac{105}{128 s^2}+\frac{315}{1024 s^3}-\frac{3968}{3 s^4 }\geq 0,
\end{equation*}
for $s\geq 40$, and so
\begin{equation}\label{emitate-bessel}
	I_2(s) \geq \frac{e^{s}}{\sqrt{2\pi s}} \left(1-\frac{15}{8 s}\right),
\end{equation}
for $s\geq 231$.

Substituting \eqref{emitate-bessel} into \eqref{defi-Gk}, and based on the following two observations:
\begin{equation*}
	\max\left\{\frac{\sqrt{2}}{\alpha_1^{\frac{3}{2}}(1)},
\frac{\sqrt{2}}{\alpha_2^{\frac{3}{2}}(1)}\right\}=\frac{3 \sqrt{42}}{49}\approx 0.396<1,
\end{equation*}
and
\[\max\left\{\frac{15}{8 \sqrt{\alpha_1(1)}},\frac{15}{8 \sqrt{\alpha_2(1)}}\right\}=\frac{15 \sqrt{21}}{56},\]
  we derive that
 for $x_k(n)\geq 152$,
\begin{align}\label{eq-S}
	G_k(n)
	&\leq\frac{\sqrt{2}}{\alpha_k^{\frac{3}{2}}(1)} \cdot \frac{144}{x_k(n)} \cdot\frac{\exp\left({-\frac{1}{2} \sqrt{\alpha_k(1)} x_k(n)}\right) }{1-\frac{15 }{8 \sqrt{\alpha_k(1) } x_k(n)}}	\nonumber\\[6pt]
	&\leq  \frac{144}{x_k(n)} \cdot \frac{\exp\left({-\frac{\sqrt{21}}{6} x_k(n)}\right)}{1- \frac{15 \sqrt{21}}{56 x_k(n)}}.		
\end{align}
Under the following observation that
\begin{align*}
	\left(1-\frac{15 \sqrt{21}}{56 x_k(n)}\right)\left(1+\frac{3}{ x_k(n)}\right)&=1+\frac{3 \left(56-5 \sqrt{21}\right)}{56 x^2_k(n)}\left(x_k(n)-\frac{15 \sqrt{21}}{56-5 \sqrt{21}}\right)\geq 1,
\end{align*}
for $x_k(n)\geq 152$,  we find that  \eqref{eq-S} can be further bounded by:
\begin{equation}\label{eq-G(x)}
	G_k(n)    \leq 144\left(\frac{1}{x_k(n)}+\frac{3}{ x^2_k(n)} \right) \exp\left({-\frac{\sqrt{21}}{6} x_k(n)}\right).
\end{equation}
We claim that  for $x_k(n)\geq 152$,
\begin{equation}\label{eq-f(x)}
 144 \exp\left({-\frac{\sqrt{21}x_k(n)}{6} }\right) \leq\frac{1}{2x_k^5(n)}.
\end{equation}
Define
\[L(s):=288s^5\exp\left({-\frac{\sqrt{21}s}{6} }\right).\]
It is evident that
\begin{align*}
L'(s)=	48 \exp\left({-\frac{\sqrt{21}s}{6} }\right) s^4 \left(-\sqrt{21} s+30\right).
\end{align*}
Since  $L'(s)\leq 0$ for $s\geq\frac{30}{  \sqrt{21}}$, we deduce that $L(s)$ is decreasing when $s \geq \frac{30}{  \sqrt{21}}$. This implies that
\[L(x_k(n))=288x^5_k(n)\exp\left({-\frac{\sqrt{21}x_k(n)}{6} }\right)\leq L(152)<1,\]
for $x_k(n)\geq 152$. So the claim is proved.

Applying \eqref{eq-f(x)} to \eqref{eq-G(x)}, we are led to
\begin{align*}
	G_k(n)\leq \left(\frac{1}{x_k(n)}+\frac{3}{ x^2_k(n)} \right) \cdot \frac{1}{2x_k^5(n)}<\frac{1}{x^6_k(n)},
\end{align*}
 for  $x_k(n)\geq 152 $. This completes the proof. \qed

The following inequality on $I_2(s)$ is also required in the proof of   Theorem \ref{thm-Delta-1}.

\begin{thm}\label{lem-B}
 For $k=1$ or $2$	and $x_k(n)\geq152$,
\begin{align}\label{eq-lem-B}
\frac{I_2^2\left(\sqrt{\alpha_k(1) } x_k(n)\right)}{I_2\left(\sqrt{\alpha_k(1) } x_k(n-1)\right)I_2\left(\sqrt{\alpha_k(1) } x_k(n+1)\right)}>1+\frac{\pi ^4 \sqrt{\alpha_k(1)}}{9 x^3_k(n)}-\frac{1100}{x_k^4(n)}.
\end{align}
\end{thm}
\proof
From \eqref{eq-BesselI1}, we see that  for  $s\geq 231$,
\begin{align}
	I_2(s)&\geq \frac{e^s}{\sqrt{2\pi s}}\left(1-\frac{15}{8s}+\frac{105}{128s^2}+\frac{315}{1024s^3}
-\frac{3968}{3s^4}\right). \label{ineqBessel-2}\\[6pt]
	I_2(s)&\leq \frac{e^s}{\sqrt{2\pi s}}\left(1-\frac{15}{8s}+\frac{105}{128s^2}
+\frac{315}{1024s^3}+\frac{3968}{3s^4}\right),\label{ineqBessel-1}
\end{align}
For convenience,  let
\[\gamma_1(k)=\frac{15 }{8 \alpha^{\frac{1}{2}}_k(1)},\quad
\gamma_2(k)=\frac{105  }{128 \alpha_k(1)},\]
\[
\gamma_3(k)=\frac{315   }{1024 \alpha^{3/2}_k(1)},\quad
\gamma_4(k)=\frac{3968  }{3 \alpha^2_k(1)},\]
and
\begin{align}\label{defi-h}
	h_k(n)&=\frac{\left(1-\frac{\gamma _1(k)}{x_k(n)}+\frac{\gamma _2(k) }{x^2_k(n)}+\frac{\gamma _3(k) }{x^3_k(n)}-\frac{\gamma _4(k)}{x_k^4(n)}\right)^2}{\left(1-\frac{\gamma _1(k)}{x_k(n-1)}+\frac{\gamma _2(k)}{x^2_k(n-1)}+\frac{\gamma _3(k)}{x^3_k(n-1)}+\frac{\gamma _4(k)}{x^4_k(n-1)}\right) }\nonumber \\[5pt]
	&\quad \quad \times \frac{1}{\left(1-\frac{\gamma _1(k)}{x_k(n+1)}+\frac{\gamma _2(k)}{x^2_k(n+1)}+\frac{\gamma _3(k)}{x^3_k(n+1)}+\frac{\gamma _4(k)}{x^4_k(n+1)}\right)}
\end{align}

Combing \eqref{ineqBessel-2} and \eqref{ineqBessel-1}, we derive that for $x_k(n)\geq 152$,
\begin{align}\label{eq-lemB-0}
	&\frac{I_2^2\left(\sqrt{\alpha_k(1) } x_k\right)}{I_2\left(\sqrt{\alpha_k(1) } x_k(n-1)\right)I_2\left(\sqrt{\alpha_k(1) } x_k(n+1)\right)}\nonumber\\[5pt]
& \qquad\geq\frac{\sqrt{x_k(n-1) x_k(n+1)}}{x_k(n)}\exp\left({\sqrt{\alpha_k(1) } (2x_k(n)-x_k(n-1)-x_k(n+1))}\right)h_k(n).
\end{align}
From the definition \eqref{defi-xk} of $x_k(n)$,  we see that
\begin{equation}\label{eq-x-Rel}
    x_k(n -1)=\sqrt{x_k^2(n)-\frac{2\pi^2}{3}},\quad x_k(n +1)=\sqrt{x_k^2(n)+\frac{2\pi^2}{3}}.
\end{equation}
This implies that
\[0<\frac{\sqrt{x_k(n-1)x_k(n+1)}}{x_k(n)}<1,\]
and so
\begin{align}\label{eq-lemB-1}
	\frac{\sqrt{x_k(n-1)x_k(n+1)}}{x_k(n)} & >
\frac{x_k(n-1)^2x_k(n+1)^2}{x_k^4(n)}=1-\frac{4\pi^4}{9x_k^4(n)}.
\end{align}

 To estimate the remaining  parts on the right-hand side of \eqref{eq-lemB-0}, we plan to establish an upper bound and a lower bound for  $x_k(n-1)$ and $x_k(n+1)$ in terms of $x_k(n)$. Observe that for $x_k(n)\geq 3$,
\begin{align*}
	&x_k(n-1)=x_k(n)-\frac{\pi^2}{3x_k(n)}-\frac{\pi^4}{18x^3_k(n)}
-\frac{\pi^6}{54x_k^5(n)}-\frac{5\pi^8}{648x_k^7(n)}+
o\left(\frac{1}{x_k^9(n)}\right),\\[5pt]	&x_k(n+1)=x_k(n)+\frac{\pi^2}{3x_k(n)}-\frac{\pi^4}{18x^3_k(n)}+\frac{\pi^6}{54x_k^5(n)}
-\frac{5\pi^8}{648x_k^7(n)}+o\left(\frac{1}{x_k^9(n)}\right),
\end{align*}
so it is readily checked that  for $x_k(n)\geq 23$,
\begin{align}
	\tilde{w}_k(n)<&x_k(n-1)<\hat{w}_k(n), \label{xkn-1}\\[5pt]
	\tilde{y}_k(n)<&x_k(n+1)<\hat{y}_k(n), \label{xkn+1}
\end{align}
where
\begin{align}\left\{
\begin{aligned}\label{wylabel}
	\tilde{w}_k(n)&=x_k(n)-\frac{\pi^2}{3x_k(n)}-\frac{\pi^4}{18x^3_k(n)}-\frac{\pi^6}{x_k^5(n)},\\
	\hat{w}_k(n)&=x_k(n)-\frac{\pi^2}{3x_k(n)}-\frac{\pi^4}{18x^3_k(n)}-\frac{\pi^6}{54x_k^5(n)},\\
\tilde{y}_k(n)&=x_k(n)+\frac{\pi^2}{3x_k(n)}-\frac{\pi^4}{18x^3_k(n)},\\
	\hat{y}_k(n)&=x_k(n)+\frac{\pi^2}{3x_k(n)}-\frac{\pi^4}{18x^3_k(n)}+\frac{\pi^6}{54x_k^5(n)}.
\end{aligned}\right.
\end{align}

Combining \eqref{xkn-1} and \eqref{xkn+1}, we deduce that  for $x_k(n)\geq 23$,
\begin{equation*}	2x_k(n)-x_k(n-1)-x_k(n+1)>2x_k(n)-\hat{w}_k(n)-\hat{y}_k(n)=\frac{\pi^4}{9x^3_k(n)}>0.
\end{equation*}
Hence
\begin{align}\label{eq-lemB-2}
&	\exp\left({\sqrt{\alpha_k(1) } (2x_k(n)-x_k(n-1)-x_k(n+1))}\right)\nonumber \\[5pt]
&\quad >1+\sqrt{\alpha_k(1) } (2x_k(n)-x_k(n-1)-x_k(n+1))\nonumber\\[5pt]
&\quad \quad =1+  \frac{\pi^4\sqrt{\alpha_k(1) }}{ 9 x^3_k(n)}.
\end{align}

We proceed to show that for $x_k(n)\geq 62$,
\begin{equation}\label{eq-h-1}
	h_k(n)\geq 1-\frac{1000}{x_k^4(n)}.
\end{equation}
Define
\[P_k(n):=x_k^2(n)x^4_k(n-1) x^4_k(n+1) \left(x_k^4(n)-\gamma _1(k)x^3_k(n)+\gamma _2(k) x^2_k(n)+\gamma _3(k) x_k(n)-\gamma _4(k)\right)^2,
\]
\begin{align*}
\tilde{Q}_k(n):=&x_k^{10}(n) \left(x^4_k(n-1)-\gamma _1(k)x^3_k(n-1)+\gamma _2(k) x^2_k(n-1)+\gamma _3(k) x_k(n-1)+\gamma _4(k)\right)\nonumber\\[5pt]
 &\times \left(x^4_k(n+1)-\gamma _1(k)x^3_k(n+1)+\gamma _2(k) x^2_k(n+1)+\gamma _3(k) x_k(n+1)+\gamma _4(k)\right).
\end{align*}
and
\begin{align*}
{Q}_k(n):=&x_k^{10}(n) \left(x^4_k(n-1)-\gamma_1(k)x^2_k(n-1)\tilde{w}_k(n)+\gamma _2(k) x^2_k(n-1)+\gamma _3(k)\hat{w}_k(n)+\gamma _4(k)\right)\nonumber\\[5pt] &\times \left(x^4_k(n+1)-\gamma _1(k)x^2_k(n+1)\tilde{y}_k(n)+\gamma _2(k) x^2_k(n+1)+\gamma _3(k)\hat{y}_k(n)+\gamma _4(k)\right).
\end{align*}
Applying \eqref{xkn-1} and \eqref{xkn+1} yields that
\begin{equation*}
\tilde{Q}_k(n)\leq {Q}_k(n).
\end{equation*}
Moreover, it can be checked that  for $x_k(n)\geq 3$, 
\[\tilde{Q}_k(n)>0.\]
Hence we derive from  \eqref{defi-h} that for $x_k(n)\geq 3$,
\begin{align*}
h_k(n)&=\frac{P_k(n) }{\tilde{Q}_k(n)} \geq \frac{P_k(n) }{{Q}_k(n)}.
\end{align*}
To prove \eqref{eq-h-1}, we next show that for $x_k(n) \geq 62$,
\begin{equation}\label{eq-h-1-tempa}
	\frac{P_k(n)}{Q_k(n)}\geq1-\frac{1000}{x_k^4(n)},
\end{equation}
which is equivalent to
\begin{equation}\label{eq-h-1-tempb}
	x_k^4(n)(P_k(n)-Q_k(n))+1000Q_k(n)\geq 0,
\end{equation}
for $x_k(n) \geq 62$.

 From the definitions of $P_k(n)$ and $Q_k(n)$, together with \eqref{eq-x-Rel} and \eqref{wylabel}, we infer that the left-hand side of \eqref{eq-h-1-tempb} is a polynomial in $x_k(n)$ with degree $18$, and so we could write
\[x_k^4(n)(P_k(n)-Q_k(n))+1000Q_k(n)=\sum_{j=0}^{18} f_j(k) x_k^j(n).\]
Clearly,
\begin{equation*}
	x_k^4(n)(P_k(n)-Q_k(n))+1000Q_k(n)\geq -\sum_{j=0}^{16}|f_j(k)| x_k^j+f_{17}(k)x_k^{17}+f_{18}(k)x_k^{18}.
\end{equation*}
Moreover,  numerical evidence indicates that for   $0 \leq  j \leq 15$ and $x_k(n)\geq 20$,
\begin{equation*}
	-|f_j(k)|x_{k}^j(n)\geq -|f_{16}(k)|x_k^{16}(n),
\end{equation*}
and
\begin{align*}
	f_{16}(k)&=-\frac{4340  }{\alpha^3_k(1)}+\frac{20625  }{4 \alpha_k(1)}-\frac{145 \pi ^4  }{96 \alpha_k(1)},\\[6pt]
f_{17}(k)&=\frac{9920  }{\sqrt{\alpha^5_k(1)}}-\frac{3750  }{\sqrt{\alpha_k(1)}}+\frac{5 \pi ^4  }{8 \sqrt{\alpha_k(1)}},\\[6pt]
f_{18}(k)&=1000-\frac{15872  }{3  \alpha^2_k(1)}.
\end{align*}
It is readily checked that for $x_k(n)\geq 62$,
\begin{equation*}
	f_{18}(k)x_k^{2}(n)+f_{17}(k)x_k^{2}(n)-17|f_{16}(k)|\geq 0.
\end{equation*}
 Assembling all these results above,  we conclude that for $x_k(n)\geq 62$,
\begin{align*}
&x_k^4(n)(P_k(n)-Q_k(n))+1000Q_k(n)\\[6pt]
&\qquad\geq (f_{18}(k)x_k^{2}(n)+f_{17}(k)x_k(n)-17|f_{16}(k)|)x_k^{16}(n) \geq 0.
\end{align*}
This proves \eqref{eq-h-1-tempb}, and so \eqref{eq-h-1} is valid.

Applying \eqref{eq-lemB-1}, \eqref{eq-lemB-2} and  \eqref{eq-h-1} to \eqref{eq-lemB-0}, we derive that for $x_k(n)\geq 152$,
\begin{align*}
&\frac{I_2^2\left(\sqrt{\alpha_k(1) } x_k(n)\right)}{I_2\left(\sqrt{\alpha_k(1) } x_k(n-1)\right)I_2\left(\sqrt{\alpha_k(1) } x_k(n+1)\right)}\\[5pt]
 &  \quad \geq \left(1-\frac{4\pi^4}{9x_k^4(n)}\right) \left(1+   \frac{\pi^4\sqrt{\alpha_k(1) } }{ 9 x^3_k(n)}\right)\left(1-\frac{1000}{x_k^4(n)}\right) \\[6pt]
	& \quad  =1+\frac{\pi ^4 \sqrt{\alpha_k(1)}}{9 x^3_k(n)}-\frac{1000+\frac{4 \pi ^4}{9}}{x_k^4(n)}-\frac{\frac{4}{81}  \pi ^8 \sqrt{\alpha_k(1)}+\frac{1000}{9} \pi ^4 \sqrt{\alpha_k(1)}}{x_k^7(n)}\\[6pt]
	&\qquad\quad +\frac{4000 \pi ^4}{9 x_k^8(n)}  +\frac{4000 \pi ^8 \sqrt{\alpha_k(1)}}{81 x_k^{11}(n)}.
\end{align*}
It's easy to check that  for $x_k(n) \geq 7$
\begin{equation*}
	\frac{100-\frac{4 \pi ^4}{9}}
	{x_k^4(n)} -\frac{\frac{4}{81}  \pi ^8 \sqrt{\alpha_k(1)}+\frac{1000}{9} \pi ^4 \sqrt{\alpha_k(1)}}{x_k^7(n)}\geq 0.
\end{equation*}
Therefore,  we arrive at
\begin{equation}
	\frac{I^2_2\left(\sqrt{\alpha_k(1) } x_k\right)}{I_2\left(\sqrt{\alpha_k(1) } x_k(n-1)\right)I_2\left(\sqrt{\alpha_k(1) } x_k(n+1)\right)}\geq  1+\frac{\pi ^4 \sqrt{\alpha_k(1)}}{9 x^3_k(n)}-\frac{1100}{x_k^4(n)}
\end{equation}
for $x_k(n) \geq 152$. This completes the proof. \qed

With Theorem \ref{lem-Delta-1} and Theorem \ref{lem-B} in hand,  we are now in a position to  give a proof of Theorem \ref{thm-Delta-1}.

\noindent {\it Proof of Theorem \ref{thm-Delta-1}.} To prove \eqref{thm-Delta-1-eq}, it  is enough to show that
\begin{equation}\label{eq-del-eq}
	\frac{\Delta_k^2(n)}{\Delta_k(n-1)\Delta_k(n+1)}\geq 1.
\end{equation}

Utilizing Theorem \ref{lem-Delta-1}, we find that  for $x_k(n)\geq 152$,
\begin{align} \label{eq-Delta-Log}
\nonumber	&\frac{\Delta_k^2(n)}{\Delta_k(n-1)\Delta_k(n+1)}\\[5pt]
&  \quad \geq \frac{M_k^2(n)\left(1-\frac{1}{x_k^6(n)}\right)^2}{M_k(n-1)M_k(n+1)\left(1+\frac{1}{x_k^6(n-1)}\right)^2\left(1+\frac{1}{x_k^6(n+1)}\right)^2} \nonumber \\[6pt]
 	&    \quad\geq \frac{x_k^2(n-1)x_k^2(n+1)}{x_k^4(n)} \cdot \frac{I^2_2\left(\sqrt{\alpha_k(1) } x_k(n)\right)}{I_2\left(\sqrt{\alpha_k(1) } x_k(n-1)\right)I_2\left(\sqrt{\alpha_k(1)} x_k(n+1)\right)}\cdot g_k(n),
\end{align}
where
\begin{align*}
	g_k(n):=\frac{\left(1-\frac{1}{x_k^6(n)}\right)^2}{\left(1+\frac{1}{x_k^6(n-1)}\right)^2\left(1+\frac{1}{x_k^6(n+1)}\right)^2}.
\end{align*}

We  claim that for $x_k(n)\geq 75$,
\begin{equation}\label{eq-g-k}
	g_k(n)\geq 1-\frac{10}{x_k^6(n)}.
\end{equation}
Invoking \eqref{eq-x-Rel}, we obtain
\begin{equation}\label{simpy-a}
\frac{x^2_k(n-1)x^2_k(n+1)}{x_k^4(n)}=\frac{(x^2_k(n)-\frac{2\pi^2}{3})
(x^2_k(n)+\frac{2\pi^2}{3})}{x_k^4(n)}=1-\frac{4\pi^4}{9x_k^4(n)},
\end{equation}
so $g_k(n)$ can be simplified as:
\begin{equation}\label{gk-a}
	g_k(n)=\frac{\left(x_k^6(n)-1\right)^2 \left(x_k^4(n)-\frac{4\pi^4}{9}\right)^3}{x_k^{12}(n) \left(\left(x_k^2(n)-\frac{2\pi^2} {3}\right)^3+1\right) \left(\left(x_k^2(n)+\frac{2\pi^2} {3} \right)^3+1\right)}.
\end{equation}
It is easy to show that  for $x_k(n)\geq 75$,
\begin{align}\label{thm1-a}
	&\left(x_k^6(n)-1\right)^2 \left(x_k^4(n)-\frac{4\pi^4}{9}\right)^3\nonumber\\[5pt]
	 & \quad \geq\left(x_k^{12}(n)-2x_k^6(n)\right)\left(x_k^8(n)-\frac{8\pi^4}{9}x^4_k(n)\right)\left(x_k^4(n)-\frac{4\pi^4}{9}\right)\nonumber\\[6pt]
	&  \quad =x_{k}^{24}(n)-\frac{4\pi^4}{3} x_k^{20}(n)-2 x_{k}^{18}(n)+\frac{32 \pi^8}{81} x_k^{16}(n)+ \frac{8\pi^4}{3} x_k^{14}(n)- \frac{64\pi^8}{81} x_k^{10}(n)\nonumber\\[6pt]
& \quad	\geq x_{k}^{24}(n)-\frac{4\pi^4}{3} x_k^{20}(n)-2 x_{k}^{18}(n),
	\end{align}
and
\begin{align}\label{thm1-b}
	&x_k^{12}(n) \left(\left(x_k^2(n)-\frac{2\pi^2} {3}\right)^3+1\right) \left(\left(x_k^2(n)+\frac{2\pi^2} {3} \right)^3+1\right)\nonumber\\[5pt]
	&  \quad =x_k^{24}(n)-\frac{4\pi ^4 }{3}  x_k^{20}(n)+2 x_k^{18}(n)+\frac{16\pi ^8}{27}  x_k^{16}(n)+\frac{8\pi ^4}{3}  x_k^{14}(n)-\left(\frac{64\pi ^{12}}{729} -1\right) x_k^{12}(n)\nonumber\\[5pt]
	&  \quad\leq x_k^{24}(n)-\frac{4\pi ^4}{3}  x_k^{20}(n)+3 x_k^{18}(n).
\end{align}
Applying \eqref{thm1-a} and \eqref{thm1-b} to \eqref{gk-a},  we derive that for $x_k(n)\geq 75$,
\begin{align*}
	g_k(n)&\geq \frac{x_{k}^{24}(n)-\frac{4\pi^4}{3} x_k^{20}(n)-2 x_{k}^{18}(n)}{x_k^{24}(n)-\frac{4\pi ^4 }{3} x_k^{20}(n)+3 x_k^{18}(n)}\\[5pt]
	&= 1-\frac{5}{x_k^{6}(n)-\frac{4\pi ^4 }{3}  x_k^{2}(n)+3}\\[5pt]
	&\geq 1-\frac{10}{x_k^6(n)},
\end{align*}
which yields \eqref{eq-g-k}, and the claim is proved.

Substituting \eqref{eq-lem-B}, \eqref{eq-g-k} and  \eqref{simpy-a} to \eqref{eq-Delta-Log}, we derive that for $x_k(n)\geq 152$,
\begin{align*}
		\frac{\Delta_k(n)^2}{\Delta_k(n-1)\Delta_k(n+1)}& \geq \left(1-\frac{4\pi^4}{9x_k^4(n)}\right)\left(1+\frac{\pi ^4 \sqrt{\alpha_k(1)}}{9 x_k^3(n)}-\frac{1100}{x_k^4(n)}\right)\left(1-\frac{10}{x_k^6(n)}\right)\\[6pt]
		&=1+\frac{\pi ^4 \sqrt{\alpha_k(1)}}{9 x_k^3(n)}-\frac{1110+\frac{4 \pi ^4}{9}}{x_k^4(n)}-\frac{10}{x_k^6(n)}-\frac{4 \pi ^8 \sqrt{\alpha_k(1)}}{81 x_k^7(n)}\\[6pt]
		&\qquad +\frac{4400 \pi ^4}{9 x_k^8(n)}-\frac{10 \pi ^4 \sqrt{\alpha_k(1)}}{9 x_k^9(n)}+\frac{11000+\frac{40 \pi ^4}{9}}{x_k^{10}(n)}\\[6pt]
		&\qquad + \frac{40 \pi ^8 \sqrt{\alpha_k(1)}}{81 x_k^{13}(n)}-\frac{44000 \pi ^4}{9 x_k^{14}(n)}.
 \end{align*}

 It is readily checked  that for $x_k(n)\geq 73$,
 \begin{align*}
 	&\frac{\pi ^4 \sqrt{\alpha_k(1)}}{9 x_k^3(n)}-\frac{1200}{x_k^4(n)}\geq 0,\\[6pt]
 	&\frac{90-\frac{4 \pi ^4}{9}}{x_k^4(n)}-\frac{10}{x_k^6(n)}-\frac{4 \pi ^8 \sqrt{\alpha_k(1)}}{81 x_k^7(n)}\geq 0,\\[6pt]
 	&\frac{4400 \pi ^4}{9 x_k^8(n)}-\frac{10 \pi ^4 \sqrt{\alpha_k(1)}}{9 x_k^9(n)}\geq 0,
 \end{align*}
 and
 \[\frac{40 \pi ^8 \sqrt{\alpha_k(1)}}{81 x_k^{13}(n)}-\frac{44000 \pi ^4}{9 x_k^{14}(n)}\geq 0.\]
 Assembling all these results, we conclude that for $x_k(n)\geq 152$ (that is, $n \geq 3512$),
\begin{equation}\label{eq-Delta-Log-2}
	\frac{\Delta^2_k(n)}{\Delta_k(n-1)\Delta_k(n+1)} \geq 1.
\end{equation}
It is routine to check that \eqref{eq-Delta-Log-2} is true for $1\leq n\leq 3512$, and hence the proof is complete. \qed

\section{Concluding remarks}
To conclude, we mention some questions and remarks for further investigation.
The  main objection of this paper is to dig into the Tur\'an inequalities for the broken $k$-diamond partition function where $k=1$ or $2$. But the numerical evidence   suggests that the main results in this paper are also valid for  all $k\geq 1$. To wit,
\begin{con}
For $k\geq 3$, $\Delta_k(n)$ is log-concave for $n\geq 1$, that is,
	\begin{equation}
		\Delta_k^2(n)\geq \Delta_k(n-1)\Delta_k(n+1).
	\end{equation}
\end{con}
More generally,  we conjectured that for $k\geq 3$ and $d\geq 1$,  the Jensen polynomial  $J_{\Delta_k}^{d,n}(X)$ associated to $\Delta_k(n)$ is hyperbolic for all sufficiently large $n$.

 As alluded to after the proof of Theorem \ref{eq-Del-k-thm} in Section 2, Sussman's formula could not be applied to derive the explicit formula for $\Delta_k(n)$ when $k\geq 3$. Therefore,   the crucial point to solve these two conjectures is to establish explicit formulas of $\Delta_k(n)$ for  $k\geq 3$.

 \vskip 0.2cm
\noindent{\bf Acknowledgment.} This work
was supported by   the National Science Foundation of China.

\end{document}